\documentclass{article}
\usepackage{amsfonts}
\usepackage{amsmath}

\newtheorem{theorem}{Theorem}

\newtheorem{proposition}[theorem]{Proposition}

\input{tcilatex}
\begin{document}

\begin{center}
\bigskip {\large \textbf{\ }}

\bigskip

{\LARGE Existence of Hamiltonian Structure in 3D}{\large \textbf{\footnote{%
To Prof. Dr. Yavuz Nutku. The existence of local bi-Hamiltonian structures
was somehow known to him while we were working on the Darboux-Halphen
system. I remember asking him, to be sure not misunderstood what we have
discussed, \textquotedblright you mean all these systems are
bi-Hamiltonian?\textquotedblright\ He said, \textquotedblright of course
they are\textquotedblright .}}}

\bigskip

\bigskip

H. G\={u}mral

Department of Mathematics, Yeditepe University

Kay\i \c{s}da\u{g}\i\ 34750 \.{I}stanbul Turkey

\bigskip hgumral@yeditepe.edu.tr

\bigskip

\bigskip
\end{center}

\bigskip

\textbf{Abstract }

In three dimensions, the construction of bi-Hamiltonian structure can be
reduced to the solutions of a Riccati equation with the arclength coordinate
of a Frenet-Serret frame being the independent variable. Explicit
integration of conserved quantities are connected with the coefficients of
Riccati equation which are elements of the third cohomology class. All
explicitly constructed examples of bi-Hamiltonian systems are exhausted when
this class along with the first one vanishes. The latter condition provides
integrating factor for explicit integration of Hamiltonian functions. For
the Darboux-Halphen system, the Godbillon-Vey invariant is shown to arise as
obstruction to integrability of integrating factor.

\bigskip

\section{\protect\LARGE Introduction}

A Poisson structure on a manifold is defined by a skew symmetric
contravariant bilinear form subjected to the Jacobi identity expressed as
the vanishing of the Schouten bracket of Poisson tensor with itself \cite%
{lich77}-\cite{mr94}. This structure having no non-degeneracy requirement
becomes the basic underlying geometry to study non-canonical Hamilton's
equations on odd dimensional manifolds as well as the Hamiltonian structures
of nonlinear evolution equations \cite{olver}-\cite{morrison}.

The first interesting case of a completely degenerate finite dimensional
Hamiltonian structure occurs in three dimensions. Many works have been
devoted to the study of three dimensional dynamical systems with primary
concern on quantization, construction of conserved quantities, Hamiltonian
structures, integrability problems and their numerical integration using
techniques from various areas such as Poisson geometry, differential
equations, Frobenius integrability theorem and theory of foliations \cite%
{nambu}-\cite{ben}.

In \cite{eh09}, we reduced the problem of constructing Hamiltonian
structures in three dimensions to the solutions of a Riccati equation in
moving coordinates of Frenet-Serret frame. All known examples of dynamical
systems having two compatible and explicit Hamiltonian structures are
exhausted by constant solution. We concluded that in three dimensions vector
fields which are not eigenvectors of the curl operator are at least locally
bi-Hamiltonian.

In this work, we shall extend the discussion to analysis of obstructions to
the construction of global Hamiltonian structures. We shall present the
Darboux-Halphen system as an example for which the Godbillon-Vey three-form
can be explicitly worked out as an obstruction to the existence of global
structure.

\section{\protect\LARGE Hamiltonian Systems in Three Dimensions}

We shall summarize the necessary ingredients of the bi-Hamiltonian formalism
in three dimensions. See \cite{nambu}-\cite{ben} for details and examples.
For $\mathbf{x=}\left\{ x^{i}\right\} =(x,y,z)\in 
%TCIMACRO{\U{211d} }%
%BeginExpansion
\mathbb{R}
%EndExpansion
^{3}$, $t\in 
%TCIMACRO{\U{211d} }%
%BeginExpansion
\mathbb{R}
%EndExpansion
$ and overdot denoting the derivative with respect to $t$, we consider the
system of autonomous differential equations%
\begin{equation}
\overset{\cdot }{\mathbf{x}}=\mathbf{v}\left( \mathbf{x}\right)  \label{e1}
\end{equation}%
associated with a three-dimensional smooth vector field $\mathbf{v}.$ Eq.(%
\ref{e1}) is said to be Hamiltonian if the right hand side can be written as 
$\mathbf{v}\left( \mathbf{x}\right) =\Omega \left( \mathbf{x}\right) \left(
dH\left( \mathbf{x}\right) \right) $ where $H\left( \mathbf{x}\right) $ is
the Hamiltonian function and $\Omega \left( \mathbf{x}\right) $ is the
Poisson bi-vector (i.e. a skew-symmetric, contravariant two-tensor)
subjected to the Jacobi identity $\left[ \Omega \left( \mathbf{x}\right)
,\Omega \left( \mathbf{x}\right) \right] =0$ defined by the Schouten bracket 
\cite{lich77}. In coordinates, if $\partial _{i}=\partial /\partial x^{i}$,
the Poisson bi-vector is $\Omega \left( \mathbf{x}\right) =\Omega
^{jk}\left( \mathbf{x}\right) \partial _{j}\wedge \partial _{k}$, with
summation over repeated indices, and the Jacobi identity reads $\Omega
^{i[j}\partial _{i}\Omega ^{kl]}=0$ where $[jkl]$ denotes the
antisymmetrization over three indices. It follows that in three dimensions
the Jacobi identity is a single scalar equation. One can exploit the vector
calculus and the differential forms in three dimensions to have a more
transparent understanding of Hamilton`s equations as well as the Jacobi
identity.\ Using the isomorphism 
\begin{equation}
J_{i}=\varepsilon _{ijk}\Omega ^{jk}\qquad i,j,k=1,2,3  \label{e4}
\end{equation}%
between skew-symmetric matrices and (pseudo)-vectors defined by the
completely antisymmetric Levi-Civita tensor $\varepsilon _{ijk}$, we can
write the Hamilton's equations and the Jacobi identity as 
\begin{equation}
\mathbf{v}=\mathbf{J}\times \nabla H\text{ \ \ \ \ \ \ \ }\mathbf{J}\cdot
\left( \nabla \times \mathbf{J}\right) =0,  \label{e5}
\end{equation}%
respectively. In this form the Jacobi identity is equivalent to the
Frobenius integrability condition $J\wedge dJ=0$ for the one form $%
J=J_{i}dx^{i}$. It is the condition for $J$ to define a foliation of
codimension one in three dimensional space \cite{hasan},\cite{arn}-\cite%
{reinhart}. A distinguished property of Poisson structures in three
dimensions is the invariance of the Jacobi identity under the multiplication
of Poisson vector $\mathbf{J}\left( \mathbf{x}\right) $ by an arbitrary but
non-zero factor. The identities

\begin{equation}
\mathbf{J}\cdot \mathbf{v}=0,\text{ \ \ \ }\nabla H\cdot \mathbf{v}=0
\label{e11}
\end{equation}%
follows directly from the Hamilton's equations in $\left( \ref{e5}\right) $.
The second equation in (\ref{e11}) is the expression for the conservation of
Hamiltonian function. A three dimensional vector $\mathbf{v}\left( \mathbf{x}%
\right) $ is said to be bi-Hamiltonian if there exist two different
compatible Hamiltonian structures \cite{olver},\cite{magri78}. In the
notation of equation $\left( \ref{e5}\right) $, this implies%
\begin{equation}
\mathbf{v}=\mathbf{J}_{1}\times \nabla H_{2}=\mathbf{J}_{2}\times \nabla
H_{1}  \label{e11a}
\end{equation}%
for the dynamical equations. The compatibility condition for $\mathbf{J}_{1}$
and $\mathbf{J}_{2}$ is defined by the Jacobi identity for the Poisson
pencil $\mathbf{J}_{1}+c\mathbf{J}_{2}$ for arbitrary constant $c.$

\section{\protect\LARGE Frenet-Serret Frame}

Let $\left( \mathbf{t},\mathbf{n},\mathbf{b}\right) $ denote the
Frenet-Serret frame associated with a differentiable curve $t\rightarrow 
\mathbf{x}(t)$ in some domain of the three dimensional space $%
%TCIMACRO{\U{211d} }%
%BeginExpansion
\mathbb{R}
%EndExpansion
^{3}$. Throughout, $\nabla =\left( \partial _{x},\partial _{y},\partial
_{z}\right) $ will denote the usual gradient operator in local Cartesian
coordinates. Given a vector field $\mathbf{v}$, the unit tangent vector $%
\mathbf{t},$ the unit normal $\mathbf{n},$ and the unit bi-normal $\mathbf{b}
$ can be constructed as%
\begin{equation}
\begin{array}{ccccc}
\mathbf{t}\left( \mathbf{x}\right) =\frac{\mathbf{v}\left( \mathbf{x}\right) 
}{\left\Vert \mathbf{v}\left( \mathbf{x}\right) \right\Vert } & \quad & 
\mathbf{n}\left( \mathbf{x}\right) =\frac{-\mathbf{t}\times \left( \nabla
\times \mathbf{t}\right) }{\left\Vert \mathbf{t}\times \left( \nabla \times 
\mathbf{t}\right) \right\Vert } & \quad & \mathbf{b}\left( \mathbf{x}\right)
=\mathbf{t}\left( \mathbf{x}\right) \times \mathbf{n}\left( \mathbf{x}\right)%
\end{array}
\label{e15}
\end{equation}%
and they form a right-handed orthonormal frame except at points $\mathbf{x}$
corresponding to equilibrium solutions and, those vector fields $\mathbf{v}$
satisfying the condition imposed by $\mathbf{t}\times \left( \nabla \times 
\mathbf{t}\right) =0.$ This condition excludes essentially the flows with
constant unit tangent and the points $\mathbf{x}$ at which the unit normal $%
\mathbf{n}$ (hence the bi-normal $\mathbf{b}$) have zeros. That is, the
cases the Frenet-Serret frame is not well-defined. To avoid this we may
assume that%
\begin{equation}
\left( \nabla \times \mathbf{t}\right) \neq \lambda \left( \mathbf{x}\right) 
\mathbf{t}  \label{e15b}
\end{equation}%
for arbitrary nonzero function $\lambda \left( \mathbf{x}\right) $.\ That
is, we exclude the dynamical systems whose unit tangent vectors are the
eigenvectors of the curl operator \cite{moses}, \cite{benn}.

We introduce the directional derivatives along the triad $\left( \mathbf{t},%
\mathbf{n},\mathbf{b}\right) $ as 
\begin{equation}
\begin{array}{ccccc}
\partial _{s}=\mathbf{t}\cdot \nabla & \quad & \partial _{n}=\mathbf{n}\cdot
\nabla & \quad & \partial _{b}=\mathbf{b}\cdot \nabla%
\end{array}
\label{e16}
\end{equation}%
so that the variables $(s,n,b)$ are the coordinates associated with the
Frenet-Serret frame based at the point $\mathbf{x}$. Assuming the Cartesian
coordinates are functions $\mathbf{x=x}(s,n.b)$ of the Frenet-Serret
coordinates one concludes easily that the Jacobian determinant is non-zero
and hence the inverse transformation $s=s(\mathbf{x}),$ \ $n=n(\mathbf{x}),$
\ $b=b(\mathbf{x})$ exists locally, that is, in a sufficiently small
neighborhood of a given point $\mathbf{x}_{0}\in 
%TCIMACRO{\U{211d} }%
%BeginExpansion
\mathbb{R}
%EndExpansion
^{3}$. These functions may be obtained by integrating the one-forms 
\begin{equation}
\tau =\mathbf{t}\cdot d\mathbf{x},\ \ \eta =\mathbf{n}\cdot d\mathbf{x},\ \
\beta =\mathbf{b}\cdot d\mathbf{x}  \label{oneforms}
\end{equation}%
the last two of which implies $n=$constant and $b=$constant when restricted
to the curve $\mathbf{x(}t\mathbf{)}$. By inverting equations $\left( \ref%
{e16}\right) $ we get the expression $\nabla =\mathbf{t}\partial _{s}+%
\mathbf{n}\partial _{n}+\mathbf{b}\partial _{b}$ for the Cartesian gradient
in the Frenet-Serret frame. We finally note that Eq.(\ref{e15}) is not the
only way to construct Frenet-Serret frame, it rather shows the existence
under the condition in Eq.(\ref{e15b}). For the example of the
Darboux-Halphen system, we will use its symmetries to construct the moving
frame.

\section{\protect\LARGE Jacobi Identity in Frenet-Serret Frame}

It follows from the identity $\mathbf{J}\cdot \mathbf{v}=0$ that the Poisson
vector $\mathbf{J}$ has no component along the unit tangent vector $\mathbf{t%
}.$ Hence, we set 
\begin{equation}
\mathbf{J}=A\mathbf{n}+B\mathbf{b}  \label{e26}
\end{equation}%
for unknown functions $A\left( \mathbf{x}\right) $ and $B\left( \mathbf{x}%
\right) $ satisfying $A^{2}+B^{2}\neq 0$. Assuming $A\neq 0$ and defining
the function $\mu =B/A$ the Jacobi identity for $\mathbf{J=}A(\mathbf{n}+\mu 
\mathbf{b)}$ reduces to the Riccati equation%
\begin{equation}
\partial _{s}\mu =\mathcal{H}_{\mathbf{n}}+\mu \mathcal{H}_{\mathbf{nb}}+\mu
^{2}\mathcal{H}_{\mathbf{b}}  \label{ricca}
\end{equation}%
in the arclenght variable $s$. Here, we define the helicity densities%
\begin{equation*}
\ \mathcal{H}_{\mathbf{n}}=\mathbf{n\cdot }\left( \nabla \times \mathbf{n}%
\right) ,\text{ \ }\mathcal{H}_{\mathbf{nb}}=\mathbf{n\cdot }\nabla \times 
\mathbf{b}+\mathbf{b\cdot }\nabla \times \mathbf{n,}\text{ \ }\mathcal{H}_{%
\mathbf{b}}=\mathbf{b\cdot }\left( \nabla \times \mathbf{b}\right)
\end{equation*}%
associated with the Frenet-Serret triad \cite{eh09}. These functions are the
coefficients of the Cartesian volume $dx\wedge dy\wedge dz$ in the
three-forms%
\begin{equation}
\ \Omega _{\mathbf{n}}=\eta \wedge d\eta ,\text{ \ }\Omega _{\mathbf{nb}%
}=\eta \wedge d\beta +\beta \wedge d\eta ,\text{ \ }\Omega _{\mathbf{b}%
}=\beta \wedge d\beta .  \label{e33}
\end{equation}%
We refer to \cite{batch}-\cite{arnold} for physical and geometric
interpretations of these quantities in the context of hydrodynamics. The
Riccati equation $\left( \ref{ricca}\right) $ is equivalent to a linear
second order equation and hence possesses two linearly independent solutions
leading to two Poisson vectors for the dynamical systems under
consideration. The Hamiltonian form of dynamical equations implies that the
Poisson vectors obtained from solutions of Riccati equation are always
compatible. Thus, we conclude that

\begin{proposition}
\label{mm}All dynamical systems in three dimensions subjected to the
condition (\ref{e15b}) posses two compatible Poisson vectors.
\end{proposition}

\section{Local Form of {\protect\LARGE Hamiltonian Structures}}

By constant solution of Eq.(\ref{ricca}) we mean a function of normal
coordinates $n$ and $b$. In particular, such solutions of the Riccati
equation include the case

\begin{equation*}
\mathcal{H}_{\mathbf{n}}\mathbf{=}0,\ \ \mathcal{H}_{\mathbf{b}}=0,\ 
\mathcal{H}_{\mathbf{nb}}\mathbf{=}0
\end{equation*}%
the first two of which are the conditions for the normal vectors $\mathbf{n}$
and $\mathbf{b}$ to satisfy the Jacobi identity or Frobenius integrability
criterion. The last equation turns out to be the compatibility condition for
the normal vectors to form a Poisson pencil and thereby to define
bi-Hamiltonian structure. Integrability criterion implies the existence of
one-forms $\xi $ and $\varsigma $ such that $d\eta =\xi \wedge \eta $ and $%
d\beta =\varsigma \wedge \beta $. From the compatibility condition we see
that $\xi =\varsigma $ and the common integrating factor $\xi $ satisfies $%
\eta \wedge d\xi =0$, $\beta \wedge d\xi =0$. Thus, we have $d\xi =g\eta
\wedge \beta $ for some arbitrary function $g$ and this closes the algebra
of one-forms. The integrating factor $\xi $ is itself integrable, that is, $%
\xi \wedge d\xi =0$ when $g=0$ or equivalently $d\xi =0$. By Poincar\'{e}
lemma, the integrating factor is the differential of a function and this
makes possible to integrate $\eta $ and $\beta $ explicitly in the given
coordinate representation. These results can be drawn from the requirement
that the one-form%
\begin{equation*}
\Gamma _{global}=\eta +\beta \varepsilon +fd\varepsilon
\end{equation*}%
with $\varepsilon $ being a parameter, is integrable $\Gamma _{global}\wedge
d\Gamma _{global}=0$ to all orders in $\varepsilon $ \cite{hasan}.
Eventually, all the known examples of bi-Hamiltonian systems in three
dimensions are included in this case. The normal coordinates $n$ and $b$
obtained from explicit integration of the one-forms $\eta $ and $\beta $
represent the global conserved quantities. They appear arbitrarily in the
Riccati equation and hence in the Poisson structures. Conversely, if we are
given a bi-Hamiltonian dynamical system of the form $\mathbf{v=}\psi \nabla
H_{1}\mathbf{\times }\nabla H_{2}$ as, for example, in \cite{nambu}-\cite%
{ben}, the representation in the Frenet-Serret frame requires $\mathbf{t=}%
\psi \parallel \mathbf{v}\parallel ^{-1}\nabla H_{1}\mathbf{\times }\nabla
H_{2}$ which implies the orthogonality of the unit tangent vector to the
gradients of Hamiltonian functions. The normal vectors $\mathbf{n}$ and $%
\mathbf{b}$ can then be identified via orthonormalization procedure applied
to the linearly independent vectors $\nabla H_{1}$ and\textbf{\ }$\nabla
H_{2}$ taking also the constraint $\parallel \nabla H_{1}\mathbf{\times }%
\nabla H_{2}\parallel =\parallel \mathbf{v}\parallel /\psi $ into account.
Thus, the normal vectors $\mathbf{n}$ and $\mathbf{b}$ defining the
bi-Hamiltonian structure corresponds, in the globally integrable case, to
the gradients of Hamiltonian functions defining Poisson vectors. Local
structure arises when we are not able to integrate the conserved functions
explicitly. Before discussing this case, we remark that, as far as the
Frenet-Serret frame is constructable, one can still use the normal vectors
to cast the dynamical equations into formal bi-Hamiltonian form. To
summarize, we have

\begin{proposition}
The manifestly bi-Hamiltonian equation $\mathbf{t=n\times b}$ is the local
form of the bi-Hamiltonian structure in three dimensions. Using an
appropriate volume three-form $\nu $, this equation can be written as%
\begin{equation}
i_{\partial _{s}}\nu =\eta \wedge \beta  \label{decomp}
\end{equation}%
where the one-forms $\eta $ and $\beta $ are defined by Eq.(\ref{oneforms}).
\end{proposition}

In the case that the solution $\mu $ of the Riccati equation is not a
constant, at least one of the three-forms in Eq.(\ref{e33}) is non-zero.
Thus, they arise as obstructions to extent globally the local bi-Hamiltonian
structure. To our knowledge, there is only one example in the literature for
which we can construct explicitly the obstruction to existence of global
Hamiltonian structure. This is the Darboux-Halphen system which was obtained
from the theory of surfaces more than a hundred years ago and resurrected in
theoretical physics around eighties.

\section{\protect\LARGE The Darboux-Halphen System}

In 1878, Darboux \cite{darb78} obtained, in his investigation of family of
second degree surfaces in $%
%TCIMACRO{\U{211d} }%
%BeginExpansion
\mathbb{R}
%EndExpansion
^{3}$, the system of equations 
\begin{equation}
{\frac{d(x+y)}{dt}}=xy\;,\;\;{\frac{d(y+z)}{dt}}=yz\;,\;\;{\frac{d(x+z)}{dt}}%
=xz.  \label{halp}
\end{equation}%
Soon after, Halphen \cite{halp},\cite{halphen} has given the time-dependent
transformations 
\begin{equation}
(t,x^{i})\mapsto ({\frac{\alpha t+\beta }{\gamma t+\delta }},2\gamma {\frac{%
\gamma t+\delta }{\alpha \delta -\gamma \beta }}+{\frac{(\gamma t+\delta
)^{2}}{\alpha \delta -\gamma \beta }}x^{i})  \label{haltr}
\end{equation}%
with $(x^{1},x^{2},x^{3})=(x,y,z)$, which leave Eq.(\ref{halp}) invariant.
He has also obtained the solution of the system. Starting from 1979, growing
attractions on the Darboux-Halphen equations (\ref{halp}) come from the
observations that various important models of theoretical physics admit
reductions to equations (\ref{halp}) or its generalizations. It appeared in
the analysis of $SO(3)$ invariant anti-self-dual Einstein metrics in
gravitation theory \cite{gibb}, in connection with the moduli space of
two-monopole problem \cite{ah85},\cite{ah}, the complex Bianchi IX
cosmological models and the reductions of self-dual Yang-Mills fields \cite%
{cac90}, the WDVV equations of topological field theories \cite{dub96} and,
renormalization group flow of WZW theory \cite{bom08}. It manifests several
interesting properties concerning integrability. Namely, the Darboux-Halphen
system is an algebraically non-integrable system \cite{ms95},\cite{valls06}
whose general solution can be expressed in terms of elliptic integrals \cite%
{halp}.

\subsection{Halphen's symmetries}

We let the Darboux-Halphen system (\ref{halp}) be represented in a local
coordinate system $\mathbf{x}\equiv (x,y,z)$ by the vector field 
\begin{equation}
v=(yz-xy-xz){\frac{\partial }{\partial x}}+(xz-xy-yz){\frac{\partial }{%
\partial y}}+(xy-xz-yz){\frac{\partial }{\partial z}}  \label{hal}
\end{equation}%
for which the corresponding differential equations are equivalent to Eq.(\ref%
{halp}) provided we perform the replacement $t\mapsto -t/2$. For notational
convenience, we also define the autonomous vector fields

\begin{equation}
u=2(x{\frac{\partial }{\partial x}}+y{\frac{\partial }{\partial y}}+z{\frac{%
\partial }{\partial z}})\;,\;\;\;w={\frac{\partial }{\partial x}}+{\frac{%
\partial }{\partial y}}+{\frac{\partial }{\partial z}}  \label{bos}
\end{equation}%
and observe that, they satisfy the Lie bracket relations%
\begin{equation*}
\lbrack u,v]=2v,\text{ \ }[u,w]=-2w,\text{ \ }[v,w]=u
\end{equation*}%
of the three-dimensional Lie algebra $\mathfrak{sl}(2)$. We call the set $%
(v,u,w)$ of vector fields having weights $(-1,0,1)$, respectively, an $%
\mathfrak{sl}(2)-$triple. The generators of the Halphen's symmetries are the
vector fields%
\begin{equation}
V_{t}={\frac{\partial }{\partial t}}\text{, \ }U_{t}=-2t{\frac{\partial }{%
\partial t}}+u\text{, \ }W_{t}=-t^{2}{\frac{\partial }{\partial t}}-w+tu%
\text{\ }  \label{tsymm}
\end{equation}%
on the time-extended space $I\times M$ $\subset 
%TCIMACRO{\U{211d} }%
%BeginExpansion
\mathbb{R}
%EndExpansion
\times 
%TCIMACRO{\U{211d} }%
%BeginExpansion
\mathbb{R}
%EndExpansion
^{3}$. These vector fields satisfy the conditions%
\begin{equation}
\lbrack \frac{\partial }{\partial t}+v,V_{t}]=0,\text{ \ }[\frac{\partial }{%
\partial t}+v,U_{t}]=-2(\frac{\partial }{\partial t}+v),\text{ \ }[\frac{%
\partial }{\partial t}+v,W_{t}]=-2t(\frac{\partial }{\partial t}+v)
\label{symm}
\end{equation}%
to be infinitesimal time-dependent geometric symmetries of the
Darboux-Halphen vector field $v$. Eq.(\ref{symm}) guarantee that if a
function $H$ is a time-dependent conserved quantity for $v$, then $X(H)$
with $X$ being any of the vectors in Eq.(\ref{tsymm}), is also a conserved
quantity. Geometric consequences of the time-dependent symmetries are
discussed in \cite{hg99}. Below, we will use these symmetries to show the
local existence of time-independent conserved quantity. This will arise as a
result of a locally integrable one-form. However, we will not be able to
integrate it explicitly because the integrating factor turns out to be not
integrable in the sense of Frobenius.

Having the $\mathfrak{sl}(2)-$triple $(v,u,w)$ intrinsically associated with
the Darboux-Halphen phase space, one can define a Frenet-Serret frame by
orthonormalizing these vector fields and proceed as described in previous
sections. Namely, one can introduce a Poisson vector on the plane of normal
vectors, obtain a Riccati equation from the Jacobi identity and investigate
its solutions. It is obvious from the explicit form of vector fields that
the coefficients of the three-forms in Eq.(\ref{e33}) will be non-vanishing.
In spite of the explicit expressions, this computational way is quite
cumbersome and will not be much useful unless one guarantees the explicit
integration of the one-forms defining the moving coordinates. That means,
one presumes existence of global bi-Hamiltonian structure. Instead, we will
first cast the Darboux-Halphen system into the form of a local Hamiltonian
structure as described by Eq.(\ref{decomp}) and see if we can integrate the
one-forms involved, that is, if we can extend globally the local structure.
For this purpose, we will proceed by exploiting the Halphen symmetries
further.

A direct computation shows that Eq.(\ref{symm}) are equivalent to 
\begin{equation}
\lbrack \frac{\partial }{\partial t}+v,V_{M}]=[\frac{\partial }{\partial t}%
+v,U_{M}]=[\frac{\partial }{\partial t}+v,W_{M}]=0  \label{charsymm}
\end{equation}%
where the time-dependent vector fields 
\begin{equation}
V_{M}=-v\;,\;\;\;U_{M}=u+2tv\;,\;\;\;W_{M}=-w+tu+t^{2}v  \label{sl2t}
\end{equation}%
on $M\subset 
%TCIMACRO{\U{211d} }%
%BeginExpansion
\mathbb{R}
%EndExpansion
^{3}$ are the unique characteristic (or evolutionary) forms \cite{olver} of $%
V_{t},U_{t},W_{t}$ along the Darboux-Halphen field $v$. One can check easily
the following result.

\begin{proposition}
Each of the sets $(v,u,w)$, $(V_{t},U_{t},W_{t})$ and $(V_{M},U_{M},W_{M})$
of vector fields form $\mathfrak{sl}(2)-$triples.
\end{proposition}

It follows from Eqs.(\ref{tsymm}) or (\ref{sl2t}) that the function $\rho $
defined by 
\begin{eqnarray*}
\rho ^{-1} &=&\det ((1,\mathbf{v),V}_{t},\mathbf{U}_{t},\mathbf{W}_{t}) \\
&=&\det ((1,\mathbf{v),(}1,0\mathbf{),(-}2t,\mathbf{u),(}-t^{2},\mathbf{-w}+t%
\mathbf{u)}) \\
&=&\det ((1,\mathbf{v),V}_{M},\mathbf{U}_{M},\mathbf{W}_{M}) \\
&=&\det ((1,\mathbf{v),-v},\mathbf{u}+2t\mathbf{v,}-\mathbf{w}+t\mathbf{u}%
+t^{2}\mathbf{v}) \\
&=&-(\mathbf{u}\times \mathbf{v)}\cdot \mathbf{w}=-4(x-y)(y-z)(z-x)
\end{eqnarray*}%
is the so called last multiplier \cite{whit}. Here, the boldface letters $%
\mathbf{u},\mathbf{v},\mathbf{w}$ are the components of $u,v,w$,
respectively.

\begin{proposition}
$\nu (\mathbf{x})=\rho dx\wedge dy\wedge dz\;$is the invariant volume
element on the Darboux-Halphen phase space.
\end{proposition}

In other words, with respect to the volume $\nu $ all three vector fields $%
u,v,w$ are divergence-free. This will be seen easily when we complete the
algebra of dual one-forms. Let $(\beta ,\alpha ,\gamma )$ be one-forms dual
to the $\mathfrak{sl}(2)-$triple $(v,u,w)$ so that they satisfy 
\begin{equation*}
i_{v}\beta =i_{u}\alpha =i_{w}\gamma =1\text{\ \ }
\end{equation*}%
with the rest of the pairings being zero. Using the components of the
representative vector fields in Eqs.(\ref{hal}) and (\ref{bos}) of the
chosen basis we find the expressions 
\begin{equation*}
\beta =\rho \mathbf{u}\times \mathbf{w}\cdot d\mathbf{x}\text{, \ \ }\alpha
=\rho \mathbf{w}\times \mathbf{v}\cdot d\mathbf{x}\text{, \ \ }\gamma =\rho 
\mathbf{v}\times \mathbf{u}\cdot d\mathbf{x}
\end{equation*}%
from which the orthonormality relations follow easily. For these dual basis
one-forms, corresponding to the Lie bracket relations of $(v,u,w),$ the $%
\mathfrak{sl}(2)-$Maurer-Cartan equations 
\begin{equation}
d\beta =-2\alpha \wedge \beta \;,\;\;\ d\alpha =\gamma \wedge \beta
\;,\;\;\;\;d\gamma =2\alpha \wedge \gamma  \label{mcsl}
\end{equation}%
can be obtained using the invariant definition 
\begin{equation*}
2d\omega (u,v)=u(i_{v}\omega )-v(i_{u}\omega )-i_{[u,v]}\omega
\end{equation*}%
of the exterior derivative on a one-form $\omega $ and the Lie algebra
relations.

\begin{proposition}
The Darboux-Halphen dynamical system can be written as 
\begin{equation*}
i_{v}\mu =\alpha \wedge \gamma =d\gamma /2
\end{equation*}
which is the form of local Hamiltonian structure.
\end{proposition}

Similar relations for the other vector fields of the $\mathfrak{sl}(2)-$%
triple are 
\begin{equation*}
i_{u}\mu =\gamma \wedge \beta =d\alpha \text{, \ }\ \ \ \ i_{w}\mu =\beta
\wedge \alpha =d\beta /2
\end{equation*}
and these relations verify also the invariance of the volume form under the
flows of $\mathfrak{sl}(2)-$triple. It also follows that 
\begin{equation*}
\beta \wedge d\beta =0,\text{ \ \ \ }\alpha \wedge d\alpha =\rho ^{-2}\nu 
\text{, \ \ \ }\gamma \wedge d\gamma =0
\end{equation*}%
so that the one-forms $\beta $ and $\gamma $ are integrable in the sense of
Frobenius. To this end, we recall the integration of one-forms by homotopy
formula. Given a one-form $\omega =\omega _{i}(x)dx^{i}$, let $%
H(x)=\int_{0}^{1}x^{i}\omega _{i}(sx)ds$. Then, $dH=\omega $ \cite{olver}.
However, for each of the one-forms $\beta $ and $\gamma $ the integrands $%
\rho \mathbf{u}\times \mathbf{w}\cdot \mathbf{x}$ and $\rho \mathbf{v}\times 
\mathbf{u}\cdot \mathbf{x}$ vanish identically and one cannot apply the
homotopy formula for their integration.

\subsection{Non-integrable integrating factors}

The algebraic structure on the solution space of the Darboux-Halphen system
can be obtained, as in the case of globally bi-Hamiltonian structures, from
a locally integrable one-form involving an arbitrary parameter, this time to
second order. More precisely, we consider the one form%
\begin{equation*}
\Gamma _{local}=\beta +(2\alpha m-dm)\varepsilon -m^{2}\gamma \varepsilon
^{2}+(m+n\varepsilon )d\varepsilon
\end{equation*}%
and require this to satisfy the local integrability condition $\Gamma
_{local}\wedge d\Gamma _{local}=0$ to all orders in the parameter $%
\varepsilon $. The resulting equations, beyond immediate consequences of the
Maurer-Cartan equations, include relations for integration of locally
integrable one-forms. Referring to \cite{hasan} for details of this
computation, we find that the one-forms 
\begin{equation}
\ \ \ \ \ \ \ e^{-2\int \alpha }\gamma \text{, }\ \ \ \ \ e^{2\int \alpha
}\beta  \label{closed}
\end{equation}%
are closed. Provided these are not elements of the first cohomology class
and, by Poincar\'{e} lemma, they must be locally exact, that is, $e^{-2\int
\alpha }\gamma =dH_{1}$ for some function $H_{1}$ and similarly for the
other. Since, $i_{v}\gamma =0$ by duality relations, this function is a
conserved quantity of the Darboux-Halphen system and, as we pointed out
earlier, the Halphen symmetries, when acted upon $H_{1}$ produces new
conserved quantities two of which will be sufficient to establish the
bi-Hamiltonian structure. For a given parametrized path the integral of the
one-form $\alpha $ can be computed locally and hence the local conserved
quantity along the given path can be obtained. This structure remains only
local because the integrating factor $\alpha $, for explicit integration of $%
H_{1}$, is non-integrable%
\begin{equation*}
\alpha \wedge d\alpha \neq 0
\end{equation*}%
and from the $\mathfrak{sl}(2)-$Maurer-Cartan equations, we see that this
condition is equivalent to the non-vanishing of the Godbillon-Vey three-form 
\cite{tondeur},\cite{reinhart}%
\begin{equation*}
\ \alpha \wedge \beta \wedge \gamma \neq 0\text{ . }
\end{equation*}

\end{document}